\documentclass[a4paper,10pt]{article}

\usepackage[utf8]{inputenc}
\usepackage{array}
\usepackage{lmodern}

\usepackage{amsmath,amssymb,amsfonts}

\usepackage[ruled]{algorithm2e}

\usepackage{graphicx}
\graphicspath{{images/} {build/}}

\usepackage{booktabs}

\usepackage{subcaption}
\usepackage{prettyref}
\newrefformat{tab}{Table~\ref{#1}}
\newrefformat{fig}{Figure~\ref{#1}}
\newrefformat{alg}{Algorithm~\ref{#1}}
\newrefformat{sec}{Section~\ref{#1}}

\usepackage[xcolor]{changebar}
\cbcolor{red}
\usepackage{authblk}
\usepackage{setspace}

\newcommand{\helmholtz}{\textsc{Helm\-holtz}}

\newcommand{\laplace}{\textsc{Laplace}}

\newcommand{\mymat}[1]{\mathbf{#1}}
\newcommand{\varlocal}[1]{{#1}_{\mathrm{L}}}
\newcommand{\varglobal}[1]{{#1}_{\mathrm{G}}}

\newcommand{\varu}{\mymat{u}}

\newcommand{\varf}{\mymat{F}}

\newcommand{\helmopsymbol}{H}
\newcommand{\helmop}{\mymat{\helmopsymbol}}
\newcommand{\condop}{\mymat{\hat{\helmopsymbol}}}
\newcommand{\condpart}{\helmop^{\mathrm{cond}}}
\newcommand{\primpart}{\helmop^{\mathrm{prim}}}
\newcommand{\localhelmop}{\varlocal{\helmop}}

\newcommand{\transmat}[1]{\mymat{\widetilde{#1}}}
\newcommand{\transop}{\transmat{\helmop}}

\newcommand{\subscriptfont}[1]{\mathrm{#1}}
\newcommand{\inner}{\subscriptfont{I}}
\newcommand{\bound}{\subscriptfont{B}}
\newcommand{\face}{\subscriptfont{F}}
\newcommand{\eigenspace}{\subscriptfont{E}}

\newcommand{\masssymbol}{M}
\newcommand{\stiffsymbol}{K}

\newcommand{\mass}{\mymat{\masssymbol}}
\newcommand{\stiff}{\mymat{\stiffsymbol}}

\newcommand{\massii}{\mass_{\mathrm{II}}}
\newcommand{\stiffii}{\stiff_{\mathrm{II}}}
\newcommand{\trans}{\mymat{S}}
\newcommand{\transii}{\mymat{S}_{\inner\inner}}

\newcommand{\identity}{\mymat{I}}
\newcommand{\eig}{\mymat{\Lambda}}
\newcommand{\diag}{\mymat{D}}

\newcommand{\tp}[3]{#1 \otimes #2 \otimes #3}
\newcommand{\tpp}[3]{\tp{\left(#1\right)}{\left(#2\right)}{\left(#3\right)}}
\newcommand{\ptp}[3]{\left(\tp{#1}{#2}{#3}\right)}

\newcommand{\uglobal}{\varglobal{\mymat{u}}}

\newcommand{\vargather}{\mymat{R}}
\newcommand{\varscatter}{\vargather^T}

\newcommand{\metriccoeff}{\mymat{d}}

\newcommand{\masselement}[1]{\masssymbol_{#1}}
\newcommand{\stiffelement}[1]{\stiffsymbol_{#1}}
\newcommand{\elementwidthelement}[1]{h_{#1}}
\newcommand{\metriccoeffelement}[1]{d_{#1}}

\newcommand{\npoint}{n_{\mathrm{I}}}
\newcommand{\nelement}{n_{\mathrm{e}}}
\newcommand{\order}[1]{\mathcal{O}\of{#1}}

\newcommand{\faceindices}{\mathcal{I}}

\newcommand{\iwest}{\mathrm{w}}
\newcommand{\ieast}{\mathrm{e}}
\newcommand{\isouth}{\mathrm{s}}
\newcommand{\inorth}{\mathrm{n}}
\newcommand{\ibottom}{\mathrm{b}}
\newcommand{\itop}{\mathrm{t}}
\newcommand{\indices}{\mathcal{I}}

\newcommand{\of}[1]{\!\left(#1\right)}

\newcommand{\tpt}{TPT}
\newcommand{\tpc}{TPC}
\newcommand{\mmc}{MMC}

\newcommand{\solveruc}{UC}
\newcommand{\solverdc}{DC}
\newcommand{\solverbc}{BC}
\newcommand{\solverbt}{BT}

\newcommand{\aspectratio}{AR_{\mathrm{max}}}

\newcommand{\eqdot}{\quad .}
\newcommand{\eqcomma}{\quad ,}

\author[a,b]{Immo Huismann\thanks{Corresponding author: Immo.Huismann@tu-dresden.de}}
\author[a,b]{Jörg Stiller}
\author[a,b]{Jochen Fröhlich}
\affil[a]{Institute of Fluid Mechanics, TU Dresden}
\affil[b]{Center for Advancing Electronics Dresden (cfaed)}
\affil[{ }]{George-Bähr-Straße 3c, 01062 Dresden, Germany}

\title{Factorizing the factorization -- a spectral-element solver for elliptic equations with linear operation count}

\providecommand{\keywords}[1]{\textbf{{Keywords:}} #1}

\begin{document}
\maketitle

\begin{abstract}
  High-order methods gain more and more attention in computational fluid dynamics.
  However, the potential advantage of these methods depends critically on the availability of efficient elliptic solvers.
  With spectral-element methods, static condensation is a common approach to reduce the number of degree of freedoms and to improve the condition of the algebraic equations.
 The resulting system is block-structured and the face-based operator well suited for matrix-matrix multiplications.
 However, a straight-forward implementation scales super-linearly with the number of unknowns and, therefore, prohibits the application to high polynomial degrees.
 This paper proposes a novel factorization technique, which yields a linear operation count of just $13 N$ multiplications, where $N$ is the total number of unknowns.
 In comparison to previous work it saves a factor larger than ${3}$ and clearly outpaces unfactored variants for all polynomial degrees.
 Using the new technique as a building block for a preconditioned conjugate gradient method resulted in a runtime scaling linearly with $N$ for polynomial degrees $2\leq p \leq 32$.
 Moreover the solver proved remarkably robust for aspect ratios up to $128$.
\end{abstract}

\section{Introduction}
\label{sec:introduction}

Many problems in computational fluid dynamics are posed on Cartesian grids.
First, there are flows in rectangular domains.
Examples include the driven cavity flow in two and three dimensions~\cite{botella_1998_spectral} as well as boundary layer, plane channel, \textsc{Couette} flows, etc.
Second, a recent trend is to use Cartesian grids together with immersed boundary methods~(IBM) or cut cell methods to represent complex geometries~\cite{lai_2000_ibm}.
Third, IBM as well as the phase-field methods, the level-set methods, or the volume-of-fluid methods, frequently employ Cartesian background grids to compute the continuous phase in a multi-phase flow, while the disperse phase moving through the computational domain is represented with the chosen method~\cite{uhlmann_2005_ibm,liu_2003_phasefield}.
And lastly, similar approaches are utilized for fluid-structure coupling~\cite{aland_2010_phasefield,peskin_2002_ibm}.

Even on Cartesian grids computational methods for incompressible flows typically spend a considerable amount of the runtime in solvers for elliptic equations~\cite{ferziger_2002_cfd}.
This is also experienced in own work~\cite{kempe_2015_gpu}.
Algorithmic benefits for these solvers directly lead to noticeable performance gains and for low-order methods, extremely fast solvers are readily available~\cite{falgout_2002_hypre}.

Spectral \textsc{Fourier} methods~\cite{gottlieb_1977_spectral} constitute the optimum efficiency for the high accuracy computation of regular solutions due to the spectral convergence and the availability of fast transformations, but require periodic boundary conditions.
Spectral methods based on more general orthogonal polynomials are also employed for flow simulation~\cite{canuto_2006_spectral}, but are restricted to regular solutions and a reduced set of boundary conditions as well.
Other high-order schemes like the discontinuous \textsc{Galerkin}~(DG) methods or the spectral-element methods~(SEM) provide more geometrical flexibility and the possibility to adjust the order of approximation.
For these (and other) reasons, the latter received vital interest from the community during the last years.
Yet, fast solvers for these methods are still a matter of research.

As with high-order methods the number of degrees of freedom inside an element scales with the polynomial degree to the power of three, ways to reduce the algebraic problem size are sought.
The static-condensation method is often used to this end.
For instance the first work on SEM~\cite{patera_1984_sem} employs it, as do more recent ones~\cite{couzy_1995_condensation,yakovlev_2015_hdg}.
Other applications of static condensation include the explicit solution for cuboidal geometries~\cite{kwan_2007_condensation}, $p$-multigrid techniques for cuboidal geometries~\cite{haupt_2013_multigrid} and the application as preconditioner for a DG scheme~\cite{hartmann_2009_dg}.

The references above benefit from the static condensation with spectacular increases in performance.
However, they all share one downside:
When increasing the polynomial degree, the operation count scales super-linearly with the number of degrees of freedom, so that the method becomes less and less efficient with higher and higher polynomial degrees.
To remain efficient at high polynomial degrees, linear complexity is required throughout the entire solver, from operator execution to preconditioner to the remaining operations inside an iteration.

As the implicit treatment of diffusion terms and pressure-velocity coupling in solvers for incompressible fluid flow often reduce to a \helmholtz\ equation, the goal of this article is to provide a \helmholtz\ solver with linear scaling.
It bases on~\cite{huismann_2015_condensation} and~\cite{huismann_2014_condensation}, where preliminary variants of the condensed \helmholtz\ operator with linear operation count were derived.
While these variants resulted in linear execution times of the iterations, they outperformed unfactorized versions implemented via dense matrix-matrix multiplications only for polynomial degrees~${p>10}$.
Current simulations, however, tend to use lower polynomial degrees~\cite{beck_2014_dg,yakovlev_2015_hdg,lombard_2015_sim} so that a gain is often not achieved.
This article proposes an efficient static condensation method for a spectral-element discretization, outperforming matrix-matrix multiplications down to a polynomial degree of~{$p=2$}.

\prettyref{sec:tensor} focuses on tensor-product matrices as a necessary prerequisite and the third section on the spectral-element method.
Section~4 recapitulates the static condensation and the fifth section operators from~\cite{huismann_2015_condensation}.
Section~6, finally, puts these elements together and proposes the new method.
In Section~7 and 8, the efficiency of the new operators and solvers is quantitatively assessed with suitable test cases.

\section{Tensor-product matrices}
\label{sec:tensor}
Many partial differential equations exhibit a separable substructure~\cite{lynch_1964_tensors}, i.e.\,the differential operator can be decomposed into smaller operators acting in single coordinate directions only.
This allows for further analysis of the operator and the resulting system matrices.
Indeed, it is the basic property to lower the operation count here, as illustrated by the following very simple example.
Assume that a two-dimensional problem is discretized using a spectral method with~$n$\ degrees of freedom in each direction such that the vector of discrete unknowns is~{$\mymat{v} \in \mathbb{R}^{n^2}$}.
The system matrix~${\mymat{C} \in \mathbb{R}^{n^2,n^2}}$ is dense, so that its straightforward application requires~{$n^4$} multiplications.
If it is a tensor-product matrix, however, its application can be reformulated as
\begin{align}
  \mymat{C} \mymat{v} &= \left(\mymat{B} \otimes \mymat{A}\right) \mymat{v} = \left(\mymat{B} \otimes \identity\right) \left(\identity \otimes \mymat{A}\right) \mymat{v}
\end{align}
with~${\identity \in \mathbb{R}^{n,n}}$\ being the identity matrix, ${\mymat{A} \in \mathbb{R}^{n,n}}$ the operator in the first direction and ${\mymat{B} \in \mathbb{R}^{n,n}}$ the operator in the second one.
The consecutive application of~$\mymat{A}$\ and~$\mymat{B}$\ then requires only~${2 n^3}$ multiplications.
In general, tensor products of dimension~$d$\ require only~${d n^{d+1}}$ multiplications compared to~${n^{2d}}$ for the application of the whole matrix.

Tensor-product matrices possess additional properties that allow for factorization techniques.
E.g., the multiplication of two tensor-product matrices is reducible to the multiplication of the respective submatrices
\begin{align}
  \left(\mymat{A} \otimes \mymat{C}\right) \left(\mymat{B}\otimes\mymat{D}\right) &= \left(\mymat{A} \mymat{B}\right) \otimes \left(\mymat{C} \mymat{D}\right) \eqdot \label{eq:tensor_product_multiplication}
\intertext{Additionally,}
  \left(\mymat{A} \otimes \mymat{B}\right)^T &= \mymat{A}^T \otimes \mymat{B}^T \label{eq:tensor_product_transpose} \\
  \left(\mymat{A} \otimes \mymat{B}\right)^{-1} &= \mymat{A}^{-1} \otimes \mymat{B}^{-1} \eqdot
\end{align}
Further properties and applications of tensor-product matrices are presented in~\cite{lynch_1964_tensors} and~\cite{deville_2002_sem}.

\section{The spectral-element method for cuboidal elements}
\label{sec:method}

This paper is concerned with the \helmholtz\ equation.
The continuous problem to solve in a domain~${\Omega}$ reads
\begin{align}
  \lambda u - \Delta u &= f \label{eq:helmholtz_equation}  \eqcomma
\end{align}
where~$u$\ is the variable to solve for, ${\lambda\geq 0}$ is a real constant parameter,~$\Delta$\ the \laplace\ operator and~$f$\ the right-hand side.
This equation was first formulated in the field of acoustic research with~${\lambda < 0}$.
Nonetheless, the case {$\lambda\geq 0$} is commonly referred to as \helmholtz\ equation in the fluid dynamics community.

Decomposing the domain into~$\nelement$\ elements~$\Omega_e$, the spectral-element method leads to the discrete equation system
\begin{align}
  \vargather \localhelmop \varscatter \uglobal &= \vargather \varlocal{\mymat{F}} \eqcomma
\end{align}
where~${\varglobal{\mymat{u}}}$ is the solution vector,~${\varlocal{\mymat{F}}}$ is the discretized right-hand side,~${\vargather}$ gathers the contributions from the elements, and its transpose~{$\varscatter$} scatters the global degrees of freedom to those local to the elements~\cite{deville_2002_sem,karniadakis_2005_sem}.
The local \helmholtz\ operator~${\localhelmop}$\ is a block-diagonal matrix consisting of the element \helmholtz\ operators~${\helmop_e}$.

This paper only considers the case of cuboidal elements.
A three-dimensional tensor-product basis is utilized in each element, allowing the standard element basis functions~${\phi_{ijk}}$ to be decomposed into three one-dimensional basis functions~\cite{deville_2002_sem} such that
\begin{align}
  \forall\quad 0 \leq i,j,k \leq p:\quad
  \phi_{ijk}\of{\xi_1,\xi_2,\xi_3} &= \varphi_i\of{\xi_1} \varphi_j\of{\xi_2} \varphi_k \of{\xi_3} \label{eq:tensor_product_basis}
\end{align}
with~${\xi_{1},\xi_{2},\xi_{3} \in \left[-1,1\right]}$ the local coordinates in the element and~${\left\{ \phi_i: i = 0,\dots,p\right\}}$ the set of basis functions employed in all three directions.
These one-dimensional functions result in the standard element mass and stiffness matrices
\begin{align}
  \masssymbol _{ij} &= \int \limits _{-1}^{1} \varphi_i        \of{\xi} \varphi_j       \of{\xi} \mathrm{d}\xi \\
  \stiffsymbol_{ij} &= \int \limits _{-1}^{1} \varphi_i^{\prime} \of{\xi} \varphi_j^{\prime}\of{\xi} \mathrm{d}\xi \eqcomma
\end{align}
where the prime denotes differentiation.
Using a tensor-product basis, the operations on each element can be decomposed and the \helmholtz\ operator becomes
\begin{align}
  \begin{aligned} \helmop_e &= \metriccoeffelement{e,0} \tp{\mass}{\mass}{\mass}
    + \metriccoeffelement{e,1} \tp{\mass}{\mass}{\stiff}\\
    &+ \metriccoeffelement{e,2} \tp{\mass}{\stiff}{\mass}
    + \metriccoeffelement{e,3} \tp{\stiff}{\mass}{\mass}  \eqdot
  \end{aligned}\label{eq:helmholtz_operator}
\end{align}
For an element~${\Omega _e}$\ of extent~${\elementwidthelement{e,i}}$, ${i \in \left\{1,2,3\right\}}$, in the three coordinate directions the coefficients~${\metriccoeff_e}$\ are determined by
\begin{align}
  \metriccoeff_e &= \frac{\elementwidthelement{e,1} \elementwidthelement{e,2} \elementwidthelement{e,3}}{8} \begin{pmatrix} \lambda & \frac{4}{\elementwidthelement{e,1}^2} & \frac{4}{\elementwidthelement{e,2}^2} & \frac{4}{\elementwidthelement{e,3}^2} \end{pmatrix}^T \eqdot \label{eq:metriccoeffelement}
\end{align}
For convenience~${\metriccoeffelement{e,0}}$ incorporates the \helmholtz\ parameter~${\lambda}$, whereas the remaining components represent metric terms.
\begin{figure}
  \centering
  \hspace*{\fill}
  \begin{subfigure}[b]{0.3\textwidth}
    \includegraphics{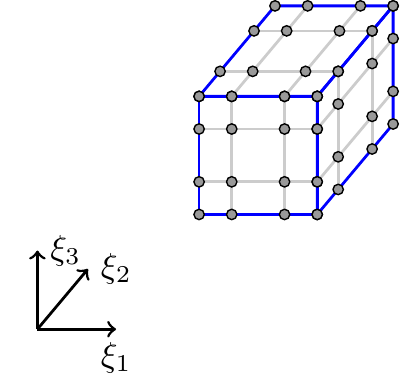}
    \caption{}
    \label{fig:tensor_product_base}
  \end{subfigure}
  \hfill
  \begin{subfigure}[b]{0.4\textwidth}
    \includegraphics{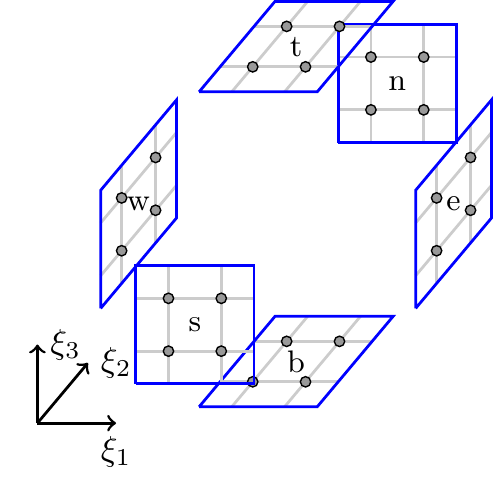}
    \caption{}
    \label{fig:compass_notation}
  \end{subfigure}
  \hspace*{\fill}
  \caption{Left: Three-dimensional tensor-product element with GLL nodes of order~${p=3}$. Right: Compass notation for element faces.}
\end{figure}%

Throughout this paper, \textsc{Lagrange} polynomials to the \textsc{Gauß-Lobatto-Legendre} (GLL) quadrature points constitute the one-dimensional nodal basis functions.
Hence, the tensor-product~GLL points can be identified with the degrees of freedom of an element, as depicted in \prettyref{fig:tensor_product_base}.
The mass matrix is approximated via GLL quadrature, generating a diagonal matrix and thereby simplifying the structure of the \helmholtz\ operator: The mass term~${\tp{\mass}{\mass}{\mass}}$\ reduces to a diagonal matrix whereas each stiffness term becomes diagonal in two dimensions, e.g.\,${\tp{\mass}{\stiff}{\mass}}$ in direction one and three.
While the simplification lowers the number of operations, it still scales with~${\order{p^4}}$\ and, hence, super-linearly with respect to the number of unknowns.
The super-linear complexity represents a roadblock on the path to higher polynomial degrees and reduces the efficiency of the solvers.
Hence, the next sections focus on the removal of this obstruction.

\section{The static condensation method}
\label{sec:condensation}

The \helmholtz\ equation is elliptic, so that the \textsc{Dirichlet} problem is well posed \cite{hirsch_1988_cfd}.
This property can be utilized to eliminate the internal nodes of an element, reducing the number of unknowns significantly, but coupling the remaining ones tighter.
As only a single element needs to be discussed, the subscript~$e$\ for element~${\Omega_e}$\ is omitted in the following.

\begin{figure}
  \hspace*{\fill}
  \begin{subfigure}{0.3\textwidth}
    \includegraphics[width=\textwidth]{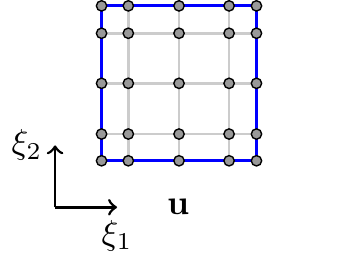}
    \caption{}
    \label{fig:categorization_full}
  \end{subfigure}
  \hfill
  \begin{subfigure}{0.3\textwidth}
    \includegraphics[width=\textwidth]{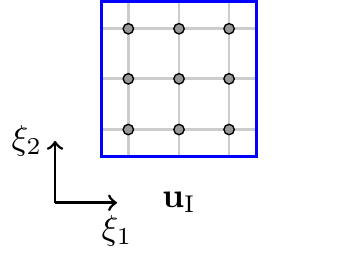}
    \caption{}
    \label{fig:categorization_inner}
  \end{subfigure}
  \hfill
  \begin{subfigure}{0.3\textwidth}
    \includegraphics[width=\textwidth]{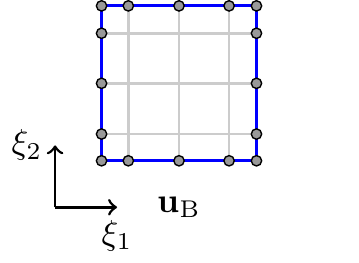}
    \caption{}
    \label{fig:categorization_bound}
  \end{subfigure}
  \hspace*{\fill}
  \caption{Categorization of the degrees of freedom in a two-dimensional element with ${p=4}$ into inner and boundary degrees. \subref{fig:categorization_full}:~all degrees of freedom in the element, \subref{fig:categorization_inner}:~nodes corresponding to internal degrees of freedom, \subref{fig:categorization_bound}:~only boundary nodes.}
  \label{fig:categorization}
\end{figure}

The degrees of freedom of the element are sorted into those located on the boundary, denoted by the subscript~$\bound$, and those in the interior, denoted by the subscript~$\inner$, as \prettyref{fig:categorization} illustrates for a two-dimensional element.
The subscripts~$\inner$\ and~$\bound$\ are also used for the corresponding matrices where applicable, e.g., ${\helmop_{\inner\bound}}$\ stands for the part of the \helmholtz\ operator mapping from the boundary to the inner degrees of freedom whereas~${\massii}$\ refers to the inner part of the one-dimensional mass matrix.%
After sorting the variables accordingly, the element \helmholtz\ operator becomes
\begin{align}
  &\begin{pmatrix}
    \helmop _{\bound \bound} & \helmop _{\bound \inner} \\
    \helmop _{\inner \bound} & \helmop _{\inner \inner}
  \end{pmatrix}
                               \begin{pmatrix}
                                 \varu _{\bound} \\ \varu_{\inner}
                               \end{pmatrix}
                             =
                               \begin{pmatrix}
                                 \varf _{\bound} \\ \varf_{\inner}
                               \end{pmatrix} \eqdot
\intertext{As~${\helmop_{\inner\inner}}$\ corresponds to the homogeneous \textsc{Dirichlet} operator and, hence, is invertible, ${\varu_{\inner}}$\ equates to}
  &\varu_{\inner} = \helmop_{\inner\inner} ^{-1} \left(\varf_{\inner} - \helmop_{\inner\bound}\varu_{\bound}\right) \label{eq:condensation_inner} \eqdot
\intertext{This, in turn, yields}
  &\underbrace{\left(\helmop_{\bound\bound} - \helmop_{\bound\inner}\helmop_{\inner\inner}^{-1}\helmop_{\inner\bound}\right)}_{\condop}\varu_{\bound} =
                                                                                                                                                       \underbrace{\varf_{\bound} - \helmop_{\bound\inner}\helmop_{\inner\inner}^{-1} \varf_{\inner}}_{\hat{\varf}} \label{eq:condensation_system} \eqdot
                                                                                                                                                       \intertext{The operator}
                                                                                                                                                       &\condop = \underbrace{\helmop_{\bound\bound}}_{\primpart} - \underbrace{\helmop_{\bound\inner}\helmop_{\inner\inner}^{-1}\helmop_{\inner\bound}}_{\condpart}
\end{align}
defines a condensed element \helmholtz\ operator with the rank equal to the number of boundary points.
It is composed of two parts:
The primary part,~${\primpart}$, is the restriction of the full \helmholtz\ operator to the boundary nodes, whereas the condensed part, ${\condpart}$, stems from the condensation process and represents the interaction of the boundary values and the internal degrees of freedom.

After solving the equation system for the element boundary values~\prettyref{eq:condensation_system}, equation~\prettyref{eq:condensation_inner} defines the interior unknowns.
\prettyref{alg:solution_process_condensed} summarizes the resulting solution process.%
\begin{algorithm}
  \caption{Solution algorithm with static condensation.}
  \label{alg:solution_process_condensed}
  Restrict to boundary nodes, condense right-hand side\\
  \ForEach{$\Omega_e$}{
    $\hat{\varf}_e \gets \varf_{\bound,e} - \helmop_{\bound\inner,e} \helmop_{\inner\inner,e}^{-1} \varf_{\inner,e}$\\
    $\hat{\varu}_e \gets \varu_{\bound,e}$
  }
  {\ }\\
  $\hat{\varu} \gets \mathrm{Solution}(\hat{\vargather}\hat{\helmop}\hat{\vargather}^T\hat{\varu} = \hat{\varf})$\\
  {\ }\\
  Regain interior degrees of freedom\\
  \ForEach{$\Omega_e$}{
    $\varu_{\inner,e} \gets \helmop_{\inner\inner,e} ^{-1} \left(\varf_{\inner,e} - \helmop_{\inner\bound,e}\varu_{\bound,e}\right)$\\
    $\varu_{e} \gets
    \begin{pmatrix}
      \varu_{\bound,e}& \varu_{\inner,e}
    \end{pmatrix}$
  }%
\end{algorithm}%

Attention to detail is required in the implementation of~${\condop}$.
A naive matrix-matrix implementation requires~${\order{p^4}}$\ multiplications and, hence, scales super-linearly with the number of degrees of freedom.
In contrast, treating the primary and condensed part separately and exploiting the tensor-product structure yields linear scaling, as shown in~\cite{huismann_2014_condensation}.

\section{Sum-factorization of the condensed equation}
\label{sec:factorization}

The condensed operator~${\condop}$\ is composed of the primary and the condensed part.
The former is the restriction of the full \helmholtz\ operator to the boundary nodes and retains the tensor-product structure.
As an example, one gets for face east, ${\face_{\ieast}}$, in compass notation (\prettyref{fig:compass_notation})
\begin{align}
  &\begin{aligned}
    \helmop_{\face_{\ieast}\face_{\ieast}} &= \metriccoeffelement{0}\masselement{00}  \massii  \otimes \massii
    + \metriccoeffelement{1}\stiffelement{00} \massii  \otimes \massii \\
    &+ \metriccoeffelement{2}\masselement{00}  \massii  \otimes \stiffii
    + \metriccoeffelement{3}\masselement{00}  \stiffii \otimes \massii         \eqcomma
  \end{aligned}\label{eq:east_east_primary_part}
\intertext{%
      with~${\metriccoeffelement{0}\dots\metriccoeffelement{3}}$ from \eqref{eq:metriccoeffelement}.
      This expression is readily evaluated in~${2 \npoint^3 + \order{\npoint^2}}$ multiplications, where the number of points per face is ${\npoint^2 = (p-1)^2}$.
      The other terms of the primary part map between opposite faces and are diagonal, e.g.
      }
      &\begin{aligned}
        \helmop_{\face_{\iwest}\face_{\ieast}} &= \metriccoeffelement{1} \stiffelement{0p}  \massii  \otimes \massii \eqdot
      \end{aligned} \label{eq:east_west_primary_part}
\end{align}
Hence, the whole primary part requires~${6\cdot2\npoint^3 + \order{\npoint^2}}$\ multiplications and scales linearly with the number of degrees of freedom.

The condensed part consists of three sub-operators:
\begin{align}
  \condpart &= \helmop_{\bound\inner} \helmop_{\inner\inner}^{-1} \helmop_{\inner\bound} \eqdot \nonumber
\end{align}
Due to the diagonal mass matrix, only faces map into the interior of the element.
A first implementation consists of precomputing the face-to-face operators and using them directly, as done in \prettyref{alg:face_to_face}, where the short-hand~${\indices = \left\{\iwest,\ieast,\isouth,\inorth,\ibottom,\itop\right\}}$ is utilized for the set of faces of the element.
Since the matrices possess~${\npoint^2 \cdot \npoint^2}$\ entries, the algorithm requires~${6\npoint^2 \cdot 6\npoint^2 = 36 \npoint^4}$ multiplications.
The scaling is super-linear with respect to the number of degrees of freedom, but an implementation can benefit from optimized libraries, e.g.\,level~3 BLAS routines~\cite{dongarra_1990_blas}, mitigating this drawback.
\begin{algorithm}
  \caption{Evaluation of the condensed operator in a direct face-to-face variant.}
  \label{alg:face_to_face}
  \ForEach{$j \in \faceindices$}{
    $\mymat{v}_{\face _j} \gets \sum  _{i \in \faceindices} \condpart_{\face_j \face_i}\mymat{u}_{\face_i}$
  }
\end{algorithm}%

Reference \cite{huismann_2015_condensation} of the present authors investigated a linearly scaling variant of the operator based on tensor products.
Linear scaling is possible when exploiting the tensor-product structure of all suboperators.
For~${\helmop_{\bound\inner}}$ and~${\helmop_{\inner\bound}}$ a tensor-product structure is induced by the restriction of \eqref{eq:helmholtz_operator} to the boundary nodes.
The inverse~${\helmop_{\inner\inner}^{-1}}$ can be expressed via the fast diagonalization method~\cite{lynch_1964_tensors}:
\begin{align}
  &\helmop_{\inner\inner}^{-1} = \left(\tp{\transii}{\transii}{\transii}\right) \diag^{-1} \left(\tp{\transii}{\transii}{\transii}\right)^T
\intertext{where}
  &\diag = \metriccoeffelement{0} \tp{\identity}{\identity}{\identity} + \metriccoeffelement{1} \tp{\identity}{\identity}{\eig}
  + \metriccoeffelement{2} \tp{\identity}{\eig}{\identity} + \metriccoeffelement{3} \tp{\eig}{\identity}{\identity}  \label{eq:hii_diagonal}
\intertext{and}
  &\transii \massii  \transii ^T = \identity ,\quad \transii \stiffii \transii ^T = \eig \eqdot     \label{eq:generalized_eigenvalues}
\end{align}
The last expression defines the one-dimensional transformation matrix $\transii$\ and the matrix of eigenvalues $\eig$\ to the generalized eigenproblem for $\stiffii$\ and $\massii$.
While the diagonal matrix~$\diag$\ is constructed by tensor-product matrices, its inverse is diagonal as well but not a tensor-product matrix.
Instead of using \prettyref{alg:face_to_face}, a sum-factorization in the inner element eigenspace is utilized for the face-to-face operators, e.g.~for face west
\begin{align}
  \mymat{v}_{\face_{\iwest}}
  &= \sum _{i \in \indices}\condpart _{\face_{\iwest}\face_{i}}\mymat{u}_{\face_i} \nonumber\\
  &= \sum _{i \in \indices}\underbrace{\helmop_{\face_{\iwest}\inner}\left(\tp{\transii}{\transii}{\transii}\right)}_{\helmop_{\face_{\iwest}\eigenspace}}\diag^{-1}\underbrace{\left(\tp{\transii^T}{\transii^T}{\transii^T} \right)\helmop_{\inner\face_{i}}}_{\helmop_{\eigenspace\face_{i}}}\mymat{u}_{\face_i} \nonumber\\
  \Leftrightarrow \mymat{v}_{\face_{\iwest}} &=\helmop_{\face_{\iwest}\eigenspace}\diag^{-1} \sum _{i \in \indices} \helmop_{\eigenspace\face_{i}}\mymat{u}_{\face_i} \eqcomma
\end{align}
where the index $\eigenspace$\ denotes the inner element eigenspace.
\prettyref{tab:suboperators} lists all the suboperators~${\helmop_{\eigenspace\face_{i}}}$.
Being products of tensor-product matrices they are tensor-product matrices themselves and are applicable in~{$3 \npoint^3$} multiplications.
First computing~${\diag ^{-1}\sum _{i \in \indices} \helmop_{\eigenspace\face_{i}}\mymat{u}_{\face_i}}$, then using it to calculate the six vectors~${\mymat{v}_{\face_{j}}}$ leads to \prettyref{alg:face_to_eigenspace}.
This algorithm evaluates the condensed part in~${37\npoint^3}$ multiplications and, thus, achieves linear complexity.
In addition, the memory requirements become linear as well, since only~${\diag^{-1}}$\ is required, whereas all the other matrices are independent of the number of elements.

While \prettyref{alg:face_to_face} uses fewer multiplications for~${p\ge 3}$, the implementation was only faster for polynomial degrees~${p>10}$~\cite{huismann_2015_condensation}.
As current high-order large-scale simulations employ polynomial degrees ranging between~$8$\ and~$12$~\cite{beck_2014_dg,merzari_2013_sim}, the possible gains actually achieved are small or even negative.
Employing further factorization for \prettyref{alg:face_to_eigenspace} enables more efficient discrete operators, as developed below.

\begin{algorithm}
  \caption{Evaluation of the condensed part that accumulates contributions in the eigenspace and then maps back to the faces.}
  \label{alg:face_to_eigenspace}
  $\mymat{\check{u}} \gets \sum  _{i \in \faceindices} \helmop_{\eigenspace\face_i}\mymat{u}_{\face_i}$\\
  $\mymat{\check{v}} \gets \diag ^{-1} \mymat{\check{u}}$\\
  \ForEach{$j \in \faceindices$}{
    $\mymat{v}_{\face _j} \gets \helmop _{\face_j\eigenspace}\mymat{\check{v}}$
  }
\end{algorithm}%

\begin{table}
  \centering
  \caption{Suboperators of the factorized condensed part of one element.}
  \begin{tabular}{clllll}
    \toprule
    $\mathrm{i}$ &\phantom{a}& $\helmop _{\face_i \eigenspace}$&\phantom{a}& $\helmop _{\eigenspace \face_i}$ \\\midrule
    $\mathrm{w}$ && $\metriccoeffelement{1} \tpp{\massii\transii}{\massii\transii}{\stiff_{0\inner}\transii}$ && $\metriccoeffelement{1} \tpp{\transii^T\massii}{\transii^T\massii}{\transii^T\stiff_{\inner 0}}$ \\
    $\mathrm{e}$ && $\metriccoeffelement{1} \tpp{\massii\transii}{\massii\transii}{\stiff_{p\inner}\transii}$ && $\metriccoeffelement{1} \tpp{\transii^T\massii}{\transii^T\massii}{\transii^T\stiff_{\inner p}}$ \\
    $\mathrm{s}$ && $\metriccoeffelement{2} \tpp{\massii\transii}{\stiff_{0\inner}\transii}{\massii\transii}$ && $\metriccoeffelement{2} \tpp{\transii^T\massii}{\transii^T\stiff_{\inner 0}}{\transii^T\massii}$ \\
    $\mathrm{n}$ && $\metriccoeffelement{2} \tpp{\massii\transii}{\stiff_{p\inner}\transii}{\massii\transii}$ && $\metriccoeffelement{2} \tpp{\transii^T\massii}{\transii^T\stiff_{\inner p}}{\transii^T\massii}$ \\
    $\mathrm{b}$ && $\metriccoeffelement{3} \tpp{\stiff_{0\inner}\transii}{\massii\transii}{\massii\transii}$ && $\metriccoeffelement{3} \tpp{\transii^T\stiff_{\inner 0}}{\transii^T\massii}{\transii^T\massii}$ \\
    $\mathrm{t}$ && $\metriccoeffelement{3} \tpp{\stiff_{p\inner}\transii}{\massii\transii}{\massii\transii}$ && $\metriccoeffelement{3} \tpp{\transii^T\stiff_{\inner p}}{\transii^T\massii}{\transii^T\massii}$ \\
    \bottomrule
  \end{tabular}
  \label{tab:suboperators}
\end{table}

\section{Factorizing the factorization}
\label{sec:further_factorization}
\prettyref{tab:suboperators} assembles the tensor-product suboperators of the condensed part used in \prettyref{alg:face_to_eigenspace}.
Two thirds of the matrix operations stem from the application of~${\massii\transii}$\ or its transpose.
Eliminating these common terms lowers the multiplication count considerably and is a key to achieving better performance.
A coordinate transformation provides the easiest approach towards this goal, as it leads to new system matrices~${\transmat{\stiff}}$\ and~${\transmat{\mass}}$, which possess favorable properties.
By applying the matrix
\begin{align}
  \trans &=
           \begin{pmatrix}
             1 & 0        & 0 \\
             0 & \transii & 0 \\
             0 & 0        & 1
           \end{pmatrix}
\end{align}
to all three directions, the element \helmholtz\ operator~${\helmop_e}$ defined in~\prettyref{eq:helmholtz_operator} is transformed to
\begin{align}
  \transmat{\helmop}_e:&=\left(\tp{\trans}{\trans}{\trans}\right) \helmop_{e} \left(\tp{\trans^T}{\trans^T}{\trans^T}\right)\\\nonumber
                       &
                         \begin{aligned}
                           &=\metriccoeffelement{e,0}\tp{\left(\trans \mass \trans^T\right)}{\left(\trans \mass \trans^T\right)}{\left(\trans \mass \trans^T\right)}\\
                           &+\metriccoeffelement{e,1}\tp{\left(\trans \mass \trans^T\right)}{\left(\trans \mass \trans^T\right)}{\left(\trans \stiff \trans^T\right)}\\
                           &+\metriccoeffelement{e,2}\tp{\left(\trans \mass \trans^T\right)}{\left(\trans \stiff \trans^T\right)}{\left(\trans \mass \trans^T\right)}\\
                           &+\metriccoeffelement{e,3}\tp{\left(\trans \stiff \trans^T\right)}{\left(\trans \mass \trans^T\right)}{\left(\trans \stiff \trans^T\right)} \eqdot
                         \end{aligned} \\
\end{align}
Defining the transformed mass matrix
\begin{align}
  \transmat{\mass} &= \trans \mass \trans^T =
                     \begin{pmatrix}
                       \mass_{00} & 0 & 0\\
                       0 & \identity & 0\\
                       0 & 0 & \mass_{pp}
                     \end{pmatrix}
\intertext{and the transformed stiffness matrix,}
  \transmat{\stiff} &= \trans \stiff \trans^T = \begin{pmatrix}
    \stiff_{00} & \stiff_{0\inner}\transii^T & \stiff_{0p} \\
    \transii \stiff_{\inner 0} & \eig & \transii \stiff_{\inner p} \\
    \stiff_{p0} & \stiff_{p\inner}\transii^T & \stiff_{pp}
  \end{pmatrix} \eqcomma
\end{align}
reduces the transformed \helmholtz\ operator to
\begin{align}
  &\begin{aligned}
    \transmat{\helmop}_e
    & =\metriccoeffelement{e,0}\tp{\transmat{\mass}}{\transmat{\mass}}{\transmat{\mass}}
    +  \metriccoeffelement{e,1}\tp{\transmat{\mass}}{\transmat{\mass}}{\transmat{\stiff}}\\
    & +\metriccoeffelement{e,2}\tp{\transmat{\mass}}{\transmat{\stiff}}{\transmat{\mass}}
    +  \metriccoeffelement{e,3}\tp{\transmat{\stiff}}{\transmat{\mass}}{\transmat{\mass}} \eqdot
  \end{aligned}
\end{align}
In the transformed system both~${\transmat{\mass}_{\inner\inner} = \identity}$ and~${\transmat{\stiff}_{\inner\inner} = \eig}$ are diagonal.
As a result, the generalized eigenvalue decomposition of~${\transmat{\stiff}_{\inner\inner}}$ with respect to~${\transmat{\mass}_{\inner\inner} = \identity}$ possesses the transformation matrix
\begin{align}
  \transmat{\trans}_{\inner\inner} &= \identity\\
  \Rightarrow \transmat{\mass}_{\inner\inner} \transmat{\trans}_{\inner\inner} &= \identity \eqdot
\end{align}
The above identities simplify the suboperators from \prettyref{tab:suboperators} to those in \prettyref{tab:transformed_suboperators}, thereby lowering the operation count.
Where the condensed part of the original operator required {$13 \npoint^3 + 24 \npoint^3$} multiplications, the condensed part of the transformed system utilizes only $13 \npoint^3$\ multiplications only.
In addition, the primary part simplifies as well.
The matrix~${\helmop_{\face_{\ieast}\face_{\ieast}}}$ in \prettyref{eq:east_east_primary_part}, e.g., becomes
\begin{align}
    \transmat{\helmop}_{\face_{\ieast}\face_{\ieast}} &= \metriccoeffelement{0}\masselement{00}  \identity  \otimes \identity
    + \metriccoeffelement{1}\stiffelement{00} \identity  \otimes \identity
    + \metriccoeffelement{2}\masselement{00}  \identity  \otimes \eig
    + \metriccoeffelement{3}\masselement{00}  \eig \otimes \identity
                                                        \label{eq:east_east_primary_part_transformed}
\end{align}
and, hence, is diagonal as well.
The primary part of the transformed system now requires~${\order{\npoint^2}}$ multiplications compared to~${12 \npoint^3 + \order{\npoint^2}}$ with the original form.
\begin{table}
  \caption{Suboperators of the factorized condensed part in the transformed system.}
  \centering
  \begin{tabular}{clllll}
    \toprule
    $\mathrm{i}$ &\phantom{a}& $\transop _{\face_i \eigenspace}$&\phantom{a}& $\transop _{\eigenspace \face_i}$ \\\midrule
    $\mathrm{w}$ &&            $\metriccoeffelement{1} \tp{\identity}{\identity}{\transmat{\stiff}_{0\inner}}$ && $\metriccoeffelement{1} \tp{\identity}{\identity}{\transmat{\stiff}_{\inner 0}}$ \\
    $\mathrm{e}$ &&            $\metriccoeffelement{1} \tp{\identity}{\identity}{\transmat{\stiff}_{p\inner}}$ && $\metriccoeffelement{1} \tp{\identity}{\identity}{\transmat{\stiff}_{\inner p}}$ \\
    $\mathrm{s}$ &&            $\metriccoeffelement{2} \tp{\identity}{\transmat{\stiff}_{0\inner}}{\identity}$ && $\metriccoeffelement{2} \tp{\identity}{\transmat{\stiff}_{\inner 0}}{\identity}$ \\
    $\mathrm{n}$ &&            $\metriccoeffelement{2} \tp{\identity}{\transmat{\stiff}_{p\inner}}{\identity}$ && $\metriccoeffelement{2} \tp{\identity}{\transmat{\stiff}_{\inner p}}{\identity}$ \\
    $\mathrm{b}$ &&            $\metriccoeffelement{3} \tp{\transmat{\stiff}_{0\inner}}{\identity}{\identity}$ && $\metriccoeffelement{3} \tp{\transmat{\stiff}_{\inner 0}}{\identity}{\identity}$ \\
    $\mathrm{t}$ &&            $\metriccoeffelement{3} \tp{\transmat{\stiff}_{p\inner}}{\identity}{\identity}$ && $\metriccoeffelement{3} \tp{\transmat{\stiff}_{\inner p}}{\identity}{\identity}$ \\
    \bottomrule
  \end{tabular}
  \label{tab:transformed_suboperators}
\end{table}%

While the operator application simplifies, the pre- and post-processing steps expand due to the transformation.
This is reflected by \prettyref{alg:solution_process_transformed}, which can be used to solve the system in its transformed variant.
\begin{algorithm}
  \caption{Solution algorithm with static condensation in transformed system.}
  \label{alg:solution_process_transformed}
  Transform, restrict to boundary nodes, condense right-hand side\\
  \ForEach{$\Omega_e$}{
    $\transmat{\varu}_e \gets  \ptp{\trans^{-1}}{\trans^{-1}}{\trans^{-1}}^T \varu_e$\\
    $\transmat{\varf}_e \gets  \ptp{\trans^{-1}}{\trans^{-1}}{\trans^{-1}}^{\phantom{T}} \varf_e$\\
    $\hat{\transmat{\varf}}_e \gets \transmat{\varf}_{\bound,e} - \transmat{\helmop}_{\bound\inner,e} \transmat{\helmop}_{\inner\inner,e}^{-1} \transmat{\varf}_{\inner,e}$\\
    $\hat{\transmat{\varu}}_e \gets \transmat{\varu}_{\bound,e}$
  }
  {\ }\\
  $\hat{\transmat{\varu}} \gets \mathrm{Solution}(\hat{\vargather}\hat{\transmat{\helmop}}\hat{\vargather}^T\hat{\transmat{\varu}} = \hat{\transmat{\varf}}$)\\
  {\ }\\
  Regain interior degrees of freedom, transform back\\
  \ForEach{$\Omega_e$}{
    $\transmat{\varu}_{\inner,e} \gets \transmat{\helmop}_{\inner\inner,e} ^{-1}\left(\transmat{\varf}_{\inner,e} - \transmat{\helmop}_{\inner\bound,e}\transmat{\varu}_{\bound,e}\right)$\\
    $\transmat{\varu}_{e} \gets
    \begin{pmatrix} \transmat{\varu}_{\bound,e} &\transmat{\varu}_{\inner,e}   \end{pmatrix}$\\
    $\varu_e \gets \ptp{\trans}{\trans}{\trans}^T\transmat{\varu}_e$
  }
\end{algorithm}%

\section{Efficiency of operators}
\label{sec:performance_operators}

The previous sections presented several variants to apply the condensed \helmholtz\ operator.
The first one realizes~\prettyref{alg:face_to_face} using a single full matrix-matrix multiplication to couple the faces of the condensed element and is hence labeled \mmc.
The matrix incorporates primary and condensed part and requires~${36\npoint^4\nelement}$ multiplications for application.
The second variant implements \prettyref{alg:face_to_eigenspace} with tensor products and is labeled \tpc.
It uses~${12\npoint^3\nelement}$ multiplications for the primary and~${37 \npoint^3\nelement}$ for the condensed part.
The tensor-product variant for the transformed system is termed~\tpt\ in the following and only requires~${13 \npoint^3\nelement}$ multiplications in total.
\prettyref{tab:operator_complexities} summarizes the multiplication count of the variants and their precomputation costs.

Multiplications do not directly translate to runtime, as loading, storing, and execution take time as well.
Hence, performance tests were conducted to directly measure the efficiency of the different operators.
The polynomial degree was varied between~${2 \leq p \leq 32}$ for a constant number of elements~${\nelement = 512}$ and a \helmholtz\ parameter~${\lambda = \pi}$.

The operators were implemented in Fortran~2008\ and compiled with the Intel Fortran compiler v.\,2015, where
a single call to DGEMM served for computing the face-to-face interaction in \mmc.
The measurements were conducted on a single core of an Intel~Xeon~E5-2690.
The setup and application of each operator was repeated~101~times and the last~100~times averaged, to remove effects from library loading.
\begin{table}
  \centering
  \caption{Complexities of leading terms of application and precomputation steps for three different variants of the condensed \helmholtz\ operator.}
  \label{tab:operator_complexities}
  \renewcommand{\arraystretch}{1.15}
  \begin{tabular}{llllllll}\toprule
    Variant&\phantom{a}& Precomputation &\phantom{a}& Primary part & \phantom{a} & Condensed part\\\midrule
    \mmc && $\order{\npoint^5\nelement}$ && $56 \npoint^2 \nelement$ && $36 \npoint^4\nelement$\\
    \tpc && $\order{\npoint^3\nelement}$ && $12 \npoint^3 \nelement$ && $37 \npoint^3\nelement$\\
    \tpt && $\order{\npoint^3\nelement}$ && $68 \npoint^2 \nelement$ && $13 \npoint^3\nelement$\\\bottomrule
  \end{tabular}
\end{table}%

\begin{figure}
  \centering
  \begin{subfigure}{0.49\textwidth}
    \includegraphics{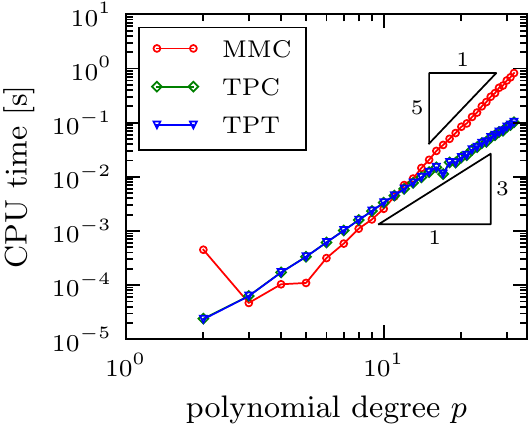}
    \caption{Precomputation}
    \label{fig:operator_efficiency_prep}
  \end{subfigure}
  \begin{subfigure}{0.49\textwidth}
    \includegraphics{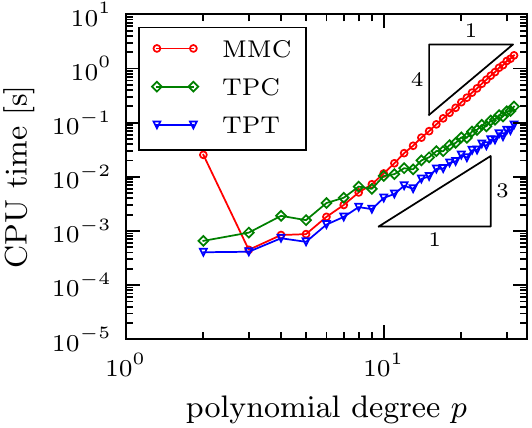}
    \caption{Operator application}
    \label{fig:operator_efficiency_comp}
  \end{subfigure}
  \caption{Average runtimes for matrix precomputation and operator evaluation for different polynomial degrees $p$\ using $\nelement = 512$. \subref{fig:operator_efficiency_prep}: Operator precomputation times, \subref{fig:operator_efficiency_comp}: operator application.}
  \label{fig:operator_efficiency}
\end{figure}%

\prettyref{fig:operator_efficiency_prep} depicts operator setup times.
The precomputation scales with the expected~${\order{\npoint^5}}$ for \mmc\ \cite{couzy_1995_condensation} and with~${\order{\npoint^3}}$ for \tpc\ and \tpt.
For homogeneous grids, where the setup of \mmc\ reduces to that of one element, the setup times are in the same order of magnitude as the time for a single application of the operator, shown in \prettyref{fig:operator_efficiency_comp}, and can therefore be neglected for unsteady solution processes requiring a huge number of time steps.
With non-uniform meshes, the setup time of \mmc\ increases~$\nelement$-fold, dominating the runtime of the whole solver process, as the setup time then is equal to the time of about~$200$ operator applications.
Furthermore, storing the full matrices becomes a problem, e.g., at double precision.
Using a polynomial degree of~${p=17}$\ and~${\nelement=512}$ the face-to-face matrices require approximately~${9.7}$ Gigabyte of memory.
The variants \tpc\ and \tpt\ do not encounter these problems due to the linearly scaling memory requirements.

\prettyref{fig:operator_efficiency_comp} depicts measured operator execution times.
For \mmc\ the operator runtime starts with some oddity at~${p=2}$, but scales with~${\npoint^4}$ starting from~${p=5}$.
\tpc\ starts with a higher runtime at~${p=3}$, but due to the lower slope it becomes faster than \mmc\ at~${p=10}$ and achieves a speedup of~${10}$ over \mmc\ at~${p=32}$.
\tpt\ exhibits the best of both worlds:
It starts with a lower runtime than \mmc\ and scales as \tpc\ does.
Furthermore, it is faster than \tpc\ by a factor of more than 2 and is~$20$\ times faster than \mmc\ at~${p=32}$.

According to \prettyref{tab:operator_complexities}, a slope of $3$\ is expected for the tensor-product based versions, but the measured runtimes exhibit a slightly smaller slope of ca.~$2.8$.
Multiple explanations are possible.
First and foremost, the primary part consists of many suboperators whose operation count scales with~${\order{\npoint^2\nelement}}$.
Second, the implementation consists of loops with~$\npoint$ iterations.
These become more efficient as the polynomial degree increases.
The combination of both can yield the lower slope.
The result is an operator whose execution time scales sub-linearly with respect to the number of degrees of freedom.

\section{Performance of solvers}
\label{sec:performance_solvers}
The previous sections focused on the linear scaling of operators as a prerequisite for solvers with linear scaling.
A typical solver, however, does not solely consist of the operator to be applied.
A good iteration scheme and preconditioner are required as well for the fast solution of the given equation.
In most cases, multigrid techniques will be employed to achieve a constant iteration count.
To investigate the impact of the condition of the system matrices on the solution procedure, a conjugate gradient solver suffices \cite{hestenes_1952_cg,shewchuk_1994_cg}, leaving only the choice of the preconditioner to be made.

As the preconditioner is required to scale linearly with respect to the number of degrees of freedom, only diagonal and tensor-product preconditioners are suited.
In \cite{huismann_2015_condensation}, block-\textsc{Jacobi} preconditioners for the faces of the elements were investigated.
These employ the exact inverse of the operators from a face to itself and can be evaluated in tensor-product form.
The remaining grid entities, i.e.\,edges and vertices, are treated similarly.
This preconditioner is referred to as the block preconditioner in the following.

Four solvers are tested here.
The first one, labeled \solveruc, is an unpreconditioned solver for the condensed system.
The second one, \solverdc, adds diagonal preconditioning.
The third one, \solverbc, applies block preconditioning to the condensed system.
All three variants utilize \tpc\ as evaluation method for the static condensed \helmholtz\ operator.
The fourth solver, labeled \solverbt, works in the transformed system and applies \tpt\ in combination with block preconditioning, which reduces to the application of a diagonal matrix in the transformed system.
\prettyref{tab:solver_complexities} summarizes the complexities of one iteration.
\begin{table}
  \centering
  \caption{Complexities of the leading terms of the different parts of the solvers investigated.}
  \label{tab:solver_complexities}
  \renewcommand{\arraystretch}{1.15}
  \begin{tabular}{lllll}\toprule
    Solver &\phantom{a}& Operator &\phantom{a}& Preconditioner\\\midrule
    \solveruc && $49\npoint^3\nelement$ && $-$ \\
    \solverdc && $49\npoint^3\nelement$ && $6 \npoint^2 \nelement$ \\
    \solverbc && $49\npoint^3\nelement$ && $24\npoint^3 \nelement$ \\
    \solverbt && $13\npoint^3\nelement$ && $6 \npoint^2 \nelement$ \\
    \bottomrule
  \end{tabular}
\end{table}

The test problem is created by a manufactured solution to \prettyref{eq:helmholtz_equation} in the domain~${\Omega = (0,2\pi)^3}$ with inhomogeneous \textsc{Dirichlet} boundary conditions.
The chosen solution is
\begin{align}
  &\begin{aligned}
    &u_{\mathrm{ex}}\of{x} = \cos\of{k (x_1 - 3 x_2 + 2 x_3)}   \sin\of{k (1 + x_1)}   \\
    &\cdot \sin\of{k (1 - x_2)}   \sin\of{k (2 x_1 + x_2)}   \sin\of{k (3 x_1 - 2 x_2 + 2 x_3)} \eqcomma
    \label{eq:solution}
  \end{aligned}
  \intertext{%
      generalizing the one employed in \cite{haupt_2013_multigrid} to three dimensions.
      The right-hand side of \prettyref{eq:helmholtz_equation} is evaluated analytically from}
    &f \of{x} = \lambda u_{\mathrm{ex}}\of{x} - \Delta u_{\mathrm{ex}}\of{x}  \eqdot
\end{align}
In the following, the case~${\lambda =  0}$ is investigated.
This in fact is the \textsc{Poisson} equation for which the resulting system matrix is not diagonally dominant anymore.
Hence, this case is harder than~${\lambda > 0}$, thus providing the ultimate test.
The stiffness parameter in~\prettyref{eq:solution} is set to ${k=5}$.

\begin{figure}
  \centering
  \hspace*{\fill}
  \begin{subfigure}{0.30\textwidth}
    \includegraphics[width=\textwidth]{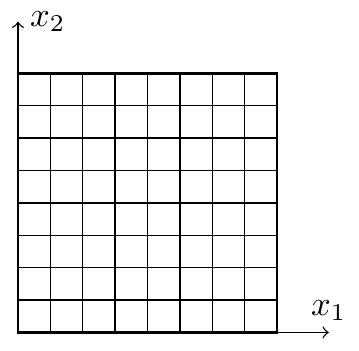}
    \caption{$\alpha = 1$}
    \label{fig:non_homogeneous_testcase_10}
  \end{subfigure}
  \hfill
  \begin{subfigure}{0.30\textwidth}
    \includegraphics[width=\textwidth]{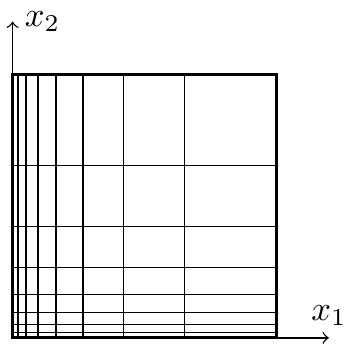}
    \caption{$\alpha = 1.5$}
    \label{fig:non_homogeneous_testcase_15}
  \end{subfigure}
  \hfill
  \begin{subfigure}{0.30\textwidth}
    \includegraphics[width=\textwidth]{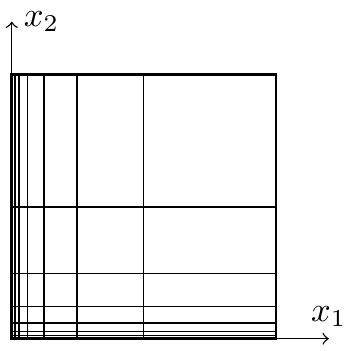}
    \caption{$\alpha = 2$}
    \label{fig:non_homogeneous_testcase_20}
  \end{subfigure}
  \hspace*{\fill}
  \caption{Meshes with~$8 \times 8$ elements in the~$x_1$-$x_2$-plane for different values of the expansion factor.}
  \label{fig:non_homogeneous_testcase}
\end{figure}

The domain is discretized using ${\nelement = 8\times 8 \times 8}$ elements, where a constant expansion factor~$\alpha$\ leads to a non-uniform spacing as illustrated in~\prettyref{fig:non_homogeneous_testcase}.
The aspect ratio of the elements in the grids can differ substantially from element to element when~${\alpha}$ gets larger.
This leads to elements of cube-like, pancake-like, and needle-like shape populating the same grid and results in a system matrix teeming with different eigenvalues due to the varying metric coefficients.
Hence, preconditioning is required to attain fast convergence and the test focuses on the gain by the preconditioner compared to the cost of applying it.

Three cases are investigated~(\prettyref{fig:non_homogeneous_testcase}): ${\alpha = 1}$, leading to a uniform mesh with a maximum aspect ratio of~${\aspectratio = 1}$,~${\alpha = 1.5}$ where the maximum aspect ratio in the~$x_1$-$x_2$~plane is~${\aspectratio \approx 17}$, and ${\alpha = 2}$ with~${\aspectratio = 128}$.
The solution process is stopped when the Euclidean norm of the residual is reduced by a factor of {$10^{12}$}.
The computations are repeated $11$ times.
Only the last $10$ repetitions contribute to the average runtime, precluding influences from initialization, e.g.\,library loading.
As for time-dependent simulations with implicit diffusion treatment the size of the time step, and thereby the \helmholtz\ parameter~$\lambda$, usually change from time step to time step if the time step size is adjusted according to a stability criterion, the precomputation times are included in the measurements.
The hardware configuration was the same as employed in \prettyref{sec:performance_operators}.

\begin{figure}
  \centering
  \begin{subfigure}{0.49\textwidth}
    \includegraphics{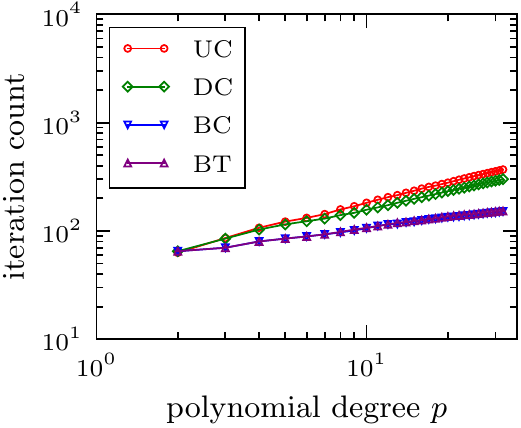}
    \caption{$\alpha = 1$}
    \label{fig:performance_non_homogeneous_po_iter_alpha_10}
  \end{subfigure}
  \begin{subfigure}{0.49\textwidth}
    \includegraphics{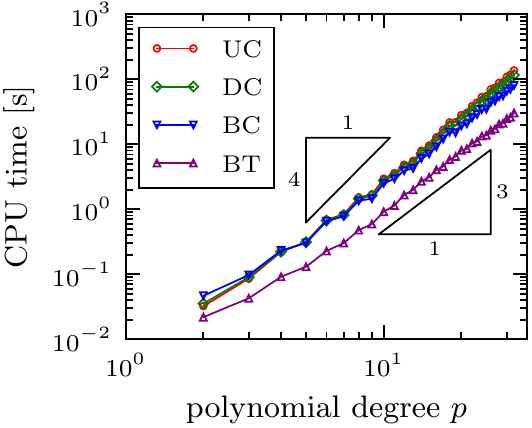}
    \caption{$\alpha = 1$}
    \label{fig:performance_non_homogeneous_po_times_alpha_10}
  \end{subfigure}

  \begin{subfigure}{0.49\textwidth}
    \includegraphics{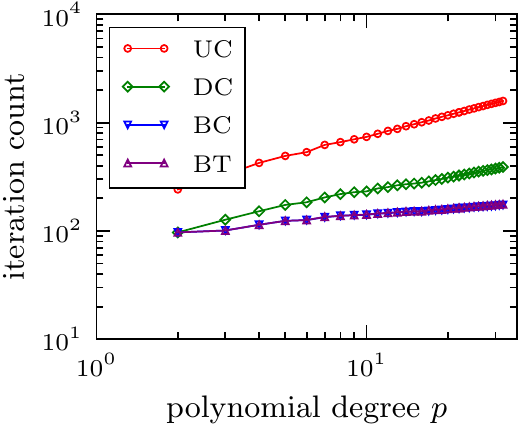}
    \caption{$\alpha = 1.5$}
    \label{fig:performance_non_homogeneous_po_iter_alpha_15}
  \end{subfigure}
  \begin{subfigure}{0.49\textwidth}
    \includegraphics{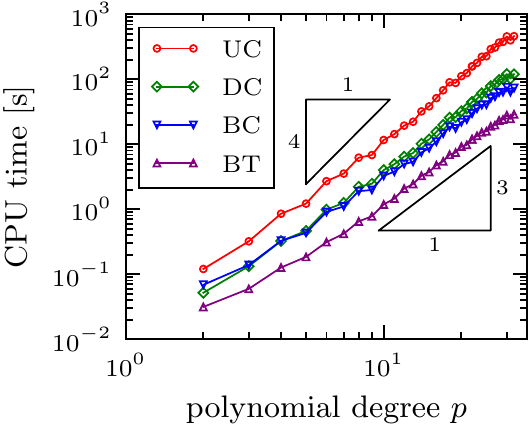}
    \caption{$\alpha = 1.5$}
    \label{fig:performance_non_homogeneous_po_times_alpha_15}
  \end{subfigure}

  \begin{subfigure}{0.49\textwidth}
    \includegraphics{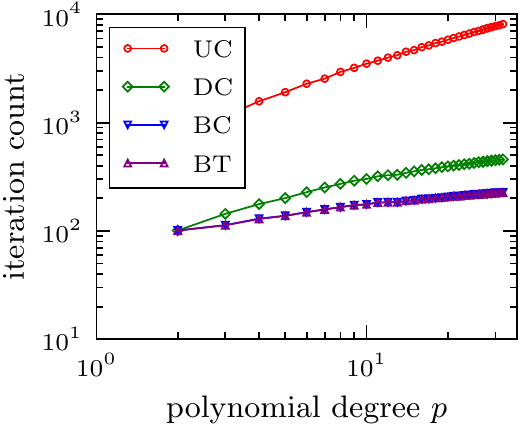}
    \caption{$\alpha = 2$}
    \label{fig:performance_non_homogeneous_po_iter_alpha_20}
  \end{subfigure}
  \begin{subfigure}{0.49\textwidth}
    \includegraphics {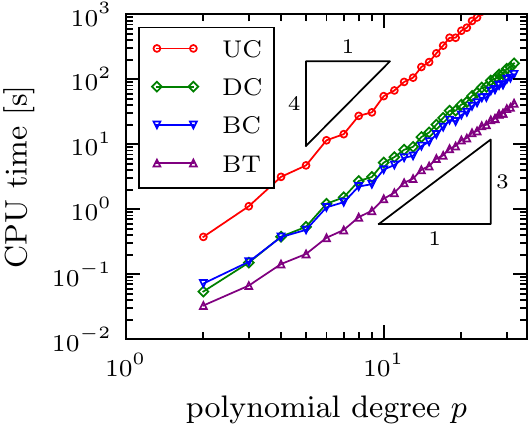}
    \caption{$\alpha = 2$}
    \label{fig:performance_non_homogeneous_po_times_alpha_20}
  \end{subfigure}

  \caption{Iteration count and CPU time for the solution of the linear system using meshes with different expansion factors.
    \subref{fig:performance_non_homogeneous_po_iter_alpha_10}~and~\subref{fig:performance_non_homogeneous_po_times_alpha_10}:~{$\alpha = 1$},
    \subref{fig:performance_non_homogeneous_po_iter_alpha_15} and \subref{fig:performance_non_homogeneous_po_times_alpha_15}:~{$\alpha = 1.5$},
    \subref{fig:performance_non_homogeneous_po_iter_alpha_20} and~\subref{fig:performance_non_homogeneous_po_times_alpha_20}:~{$\alpha = 2$}.
  }
  \label{fig:performance_non_homogeneous_po}
\end{figure}

\prettyref{fig:performance_non_homogeneous_po} summarizes the results of the test.
In all three cases, the iteration count behaves similarly:
The number of iterations starts at a low value and rises with the polynomial degree, as is to be expected.
The slope is the largest for \solveruc, slightly lower for \solverdc\ and lowest for \solverbc\ and \solverbt, which exhibit nearly the same iteration count.
For the latter three solvers, the iteration count does not differ substantially for different values of~$\alpha$, only an increase by a factor of about~$1.5$ is observed between~${\aspectratio=1}$ and~${\aspectratio=128}$.
The unpreconditioned solver is not as robust: A factor of four lies between the iteration count for~${\alpha=1}$\ and the one for~${\alpha=1.5}$ and a twenty-fold increase is found for~${\alpha=2}$.
Hence, when regarding the number of iterations, all preconditioned variants are usable, though \solverdc\ with some drawbacks compared to \solverbc\ and \solverbt.
But the unpreconditioned one is not usable for non-uniform meshes.

\prettyref{tab:iteration_times} lists the runtimes per iteration for the case of~${\alpha=2}$.
As expected, the unpreconditioned solver possesses the lowest runtime per iteration for the condensed system, with the diagonal preconditioner slightly increasing the runtime and the block preconditioner taking~${50\ \%}$\ longer per iteration.
However, this advantage is over-compensated by the iteration count of \solveruc:
The diagonal preconditioning reduces the runtime of the solver by a factor around 10 for~${p=8}$ and approximately~${17}$ for~${p=32}$.
Using the block-preconditioning yields further savings in runtime.
Yet the large effect the block-preconditioning on the iteration count is mitigated by the runtime spent for preconditioning:
The solver \solverdc\ uses only a quarter more runtime for polynomial degrees~${p\leq16}$ than \solverbc\ and requires less implementation effort.
\begin{table}
  \centering
  \caption{Computation times per degree of freedom and approximated iteration times per degree of freedom as a function of the polynomial degree obtained with the four different solvers using~${\alpha=2}$ and~${\nelement = 8^3}$.}
  \label{tab:iteration_times}
  \begin{tabular}{rrrrrrrrrrrrrrrrr}
\toprule
&&\multicolumn{7}{c}{Iteration time per DOF \mbox{$[\mathrm{ns}]$}}&&\multicolumn{7}{c}{Solution time per DOF \mbox{$[\mathrm{\mu s}]$}}\\
\cmidrule{3-9} \cmidrule{11-17} $p$ && UC && DC && BC && BT && UC && DC && BC && BT\\
\midrule
8 && 52.3 && 56.2 && 75.7 && 25.6 && 154 && 15.3 && 12.6 && 4.25 \\
12 && 33.4 && 36.8 && 49.2 && 20.3 && 133 && 12.0 && 9.05 && 3.73 \\
16 && 29.1 && 31.7 && 41.5 && 17.7 && 145 && 11.6 && 8.17 && 3.47 \\
20 && 27.3 && 29.5 && 38.9 && 16.1 && 159 && 11.6 && 8.02 && 3.31 \\
24 && 26.6 && 28.6 && 38.8 && 14.4 && 177 && 12.0 && 8.34 && 3.09 \\
28 && 24.9 && 26.5 && 36.3 && 13.1 && 185 && 11.7 && 8.03 && 2.88 \\
28 && 24.9 && 26.6 && 36.5 && 13.1 && 185 && 11.7 && 8.06 && 2.89 \\
32 && 23.6 && 25.2 && 34.4 && 12.4 && 192 && 11.5 && 7.81 && 2.81 \\
\bottomrule
\end{tabular}
\end{table}%

The solver \solverbc\ results in a large operator runtime and possesses a too costly preconditioner.
These drawbacks are removed with \solverbt:
The operator is faster and the preconditioner is diagonal in the transformed system and, hence, cheap to apply, while generating the same iteration count.
Combining both properties leads to a performance gain by a factor of $3$\ to $4$\ compared to the diagonally preconditioned case and of~$2$\ to~$3$\ compared to the block-preconditioned version.

While these savings seem insignificant, one has to keep in mind, that the solver does not solely consist of operator and preconditioner.
Many array operations are present in a CG solver and the gather-scatter operation requires runtime as well.
The new variant spends most of the time not in applying the operator nor in the preconditioner, but rather multiplying arrays etc., where no performance gain is possible and, hence, a hard barrier is present.
This also limits the potential of further factorizations.

\begin{figure}
  \centering
  \begin{subfigure}{0.49\textwidth}
    \includegraphics{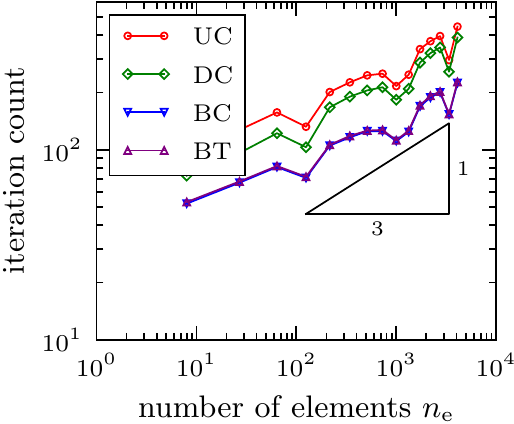}
    \caption{}
    \label{fig:performance_elements_iter}
  \end{subfigure}
  \begin{subfigure}{0.49\textwidth}
    \includegraphics{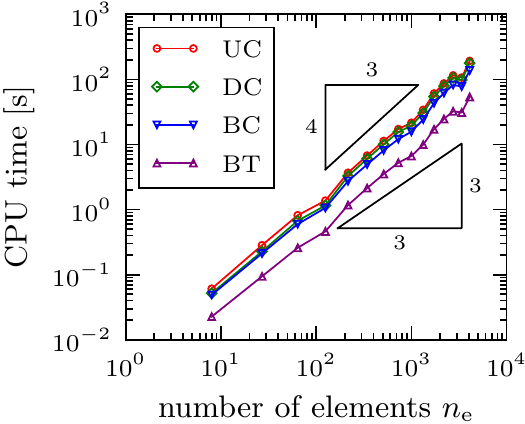}
    \caption{}
    \label{fig:performance_elements_times}
  \end{subfigure}
  \caption{Iteration count and runtimes for~${p = 16}$ when varying the number of elements.}
  \label{fig:performance_elements}
\end{figure}

To investigate the robustness of the solvers against the number of elements~${\nelement}$, the testcase~{$\alpha = 1$} was repeated for a constant polynomial degree of~${p=16}$ with~${\nelement}$ varying from~${2^3}$ to~${16^3}$.
When utilizing CG solvers, the runtime of a three-dimensional finite element solver generally scales with~${\nelement^{4/3}}$~\cite{shewchuk_1994_cg}.
\prettyref{fig:performance_elements} shows the iteration count and the CPU time of the present SEM solvers.
The number of iterations exhibits a lower slope than~${1/3}$, hence, the CPU time scales better than the expected~${\nelement^{4/3}}$.
The effect is welcomed, but the reason is probably not using enough elements to compute in the asymptotic regime.
The main conclusion from the data lies in the fact that the solvers are not robust against an increase in the number of elements.
This is to be expected, as preconditioning with the topology in mind is required to achieve this feat, e.g.\ with low-order finite elements~\cite{manna_2004_preconditioning,hartmann_2009_dg} or even multigrid~\cite{trottenberg_2001_multigrid}.

\section{Conclusions}
\label{sec:conclusions}%
This paper proposes a new evaluation technique for the condensed \helmholtz\ operator based on a tensor-product factorization, where a transformation streamlines the multiplication count.
The resulting operator variant not only scales linearly with the number of degrees of freedom, which was one major goal of this paper, but also reduces the multiplication count to a quarter of that required without the improvement~\cite{huismann_2015_condensation}.
This allows the new technique to outpace variants utilizing highly optimized libraries for matrix-matrix multiplications.
An example was provided with the \mmc\ variant that uses DGEMM.
Not only does the new method yield a speedup for all polynomial degrees over \mmc, e.g.\,by a factor of~$20$\ for~${p=32}$, it also achieves a speedup of $2$\ over previous tensor-product implementations in~\cite{huismann_2015_condensation}.

After comparing the efficiency of the operators, different solvers based on the two fastest evaluation techniques were investigated using preconditioning with linear scaling in the operation count.
Block-\textsc{Jacobi} type preconditioners provided a lower iteration count than diagonal preconditioning, which is nearly independent of the aspect ratio of the elements. For instance an increase of the maximum aspect ratio from~${\aspectratio = 1}$ to~${\aspectratio = 128}$ adds just~${50\ \%}$ to the iteration count.
Yet in the original condensed system, the block preconditioning is far more expensive than in the transformed system.
This makes the new solver $2$\ to $3$\ times faster than the previous variants.
Moreover, the scaling of the operators leads to a linearly scaling runtime of the solver with respect to the number of degrees of freedom.
This was indeed achieved with standard programming language, compiler, and hardware.

While the proposed solver scales very well with respect to the polynomial degree, the performance degrades with the number of elements.
This issue was mostly disregarded here, except with some timing measurements, because it can be removed by a multigrid approach~\cite{trottenberg_2001_multigrid}.
This is beyond the scope of this paper, which deals with the discrete representation of the \helmholtz\ operator itself.
Let us just mention that in the transformed system efficient preconditioning can be diagonal which is well suited for operator-based multigrid variants such as cascadic multigrid~\cite{deuflhard_1996_cascadic_multigrid} or multigrid CG methods~\cite{pflaum_2008_multigrid}.
Future work will focus on combining the operator factorization laid out in this paper with multigrid techniques to attain a constant iteration count and, hence, constructing an entirely linearly scaling solver.

\subsubsection*{Acknowledgment of Funds:}
This work is supported in part by the German Research Foundation (DFG) within the Cluster of Excellence ‘Center for Advancing Electronics Dresden’ (cfaed).

\bibliographystyle{abbrv}
\bibliography{transformed_condensation}

\begin{thebibliography}{10}

\bibitem{aland_2010_phasefield}
S.~Aland, J.~S. Lowengrub, and A.~Voigt.
\newblock Two-phase flow in complex geometries: A diffuse domain approach.
\newblock {\em CMES}, 57(1):77--106, 2010.

\bibitem{beck_2014_dg}
A.~D. Beck, T.~Bolemann, D.~Flad, H.~Frank, G.~J. Gassner, F.~Hindenlang, and
  C.-D. Munz.
\newblock High-order discontinuous {Galerkin} spectral element methods for
  transitional and turbulent flow simulations.
\newblock {\em International Journal for Numerical Methods in Fluids},
  76(8):522--548, 2014.

\bibitem{deuflhard_1996_cascadic_multigrid}
F.~A. Bornemann and P.~Deuflhard.
\newblock The cascadic multigrid method for elliptic problems.
\newblock {\em Numerische Mathematik}, 75(2):135--152, 1996.

\bibitem{botella_1998_spectral}
O.~Botella and R.~Peyret.
\newblock Benchmark spectral results on the lid-driven cavity flow.
\newblock {\em Computers \& Fluids}, pages 421--433, 1998.

\bibitem{canuto_2006_spectral}
C.~Canuto, M.~Hussaini, A.~Quarteroni, and T.~A. Zang.
\newblock {\em Spectral Methods. Fundamentals in Single Domains}.
\newblock Springer-Verlag Berlin Heidelberg, 2006.

\bibitem{couzy_1995_condensation}
W.~Couzy and M.~Deville.
\newblock A fast {Schur} complement method for the spectral element
  discretization of the incompressible {Navier}-{Stokes} equations.
\newblock {\em Journal of Computational Physics}, 116(1):135 -- 142, 1995.

\bibitem{deville_2002_sem}
M.~Deville, P.~Fischer, and E.~Mund.
\newblock {\em High-Order Methods for Incompressible Fluid Flow}.
\newblock Cambridge University Press, 2002.

\bibitem{dongarra_1990_blas}
J.~Dongarra, J.~Ducroz, S.~Hammarling, and I.~Duff.
\newblock {A set of level 3 basic linear algebra subprograms}.
\newblock {\em {ACM Transactions on Mathematical Software}}, {16}({1}):{1--17},
  {Mar} {1990}.

\bibitem{falgout_2002_hypre}
R.~Falgout and U.~Yang.
\newblock hypre: A library of high performance preconditioners.
\newblock In P.~Sloot, C.~Tan, J.~Dongarra, and A.~Hoekstra, editors, {\em
  Computational Science - ICCS 2002, Proceedings, Part III}, volume 2331 of
  {\em Lecture Notes in Computer Science}, pages 632--641. Springer-Verlag
  Berlin Heidelberg, 2002.
\newblock International Conference on Computational Science, Amsterdam,
  Netherlands, Apr 21-24, 2002.

\bibitem{ferziger_2002_cfd}
J.~H. Ferziger and M.~Peri\'c.
\newblock {\em Computational Methods for Fluid Dynamics}.
\newblock Springer-Verlag, 2002.

\bibitem{gottlieb_1977_spectral}
D.~Gottlieb and S.~Orszag.
\newblock {\em Numerical Analysis of Spectral Methods: Theory and
  Applications}.
\newblock CBMS-NSF Regional Conference Series in Applied Mathematics. Society
  for Industrial and Applied Mathematics, 1977.

\bibitem{hartmann_2009_dg}
R.~Hartmann, M.~Lukáčová-Medvid'ová, and F.~Prill.
\newblock Efficient preconditioning for the discontinuous {Galerkin} finite
  element method by low-order elements.
\newblock {\em Applied Numerical Mathematics}, 59(8):1737 -- 1753, 2009.

\bibitem{haupt_2013_multigrid}
L.~Haupt, J.~Stiller, and W.~E. Nagel.
\newblock A fast spectral element solver combining static condensation and
  multigrid techniques.
\newblock {\em Journal of Computational Physics}, 255(0):384 -- 395, 2013.

\bibitem{hestenes_1952_cg}
M.~R. Hestenes and E.~Stiefel.
\newblock Methods of conjugate gradients for solving linear systems.
\newblock {\em Journal of Research of the National Bureau of Standards},
  {49}({6}):{409--436}, {1952}.

\bibitem{hirsch_1988_cfd}
C.~Hirsch.
\newblock {\em Numerical Computation of Internal \& External Flows:
  Fundamentals of Numerical Discretization}.
\newblock John Wiley \& Sons, Inc., New York, NY, USA, 1988.

\bibitem{huismann_2014_condensation}
I.~Huismann, L.~Haupt, J.~Stiller, and J.~Fröhlich.
\newblock Sum factorization of the static condensed helmholtz equation in a
  three-dimensional spectral element discretization.
\newblock {\em PAMM}, 14(1):969--970, 2014.

\bibitem{huismann_2015_condensation}
I.~Huismann, J.~Stiller, and J.~Fröhlich.
\newblock Fast static condensation for the {Helmholtz} equation in a
  spectral-element discretization.
\newblock In {\em Proceedings of the 11th international conference on Parallel
  Processing and Applied Mathematics}, 2015.
\newblock To appear.

\bibitem{karniadakis_2005_sem}
G.~Karniadakis and S.~Sherwin.
\newblock {\em Spectral/hp Element Methods for Computational Fluid Dynamics}.
\newblock Oxford University Press, 2nd edition, 2005.

\bibitem{kempe_2015_gpu}
T.~Kempe, A.~Aguilera, W.~Nagel, and J.~Fröhlich.
\newblock Performance of a projection method for incompressible flows on
  heterogeneous hardware.
\newblock {\em Computers \& Fluids}, 121:37 -- 43, 2015.

\bibitem{kwan_2007_condensation}
Y.-Y. Kwan and J.~Shen.
\newblock An efficient direct parallel spectral-element solver for separable
  elliptic problems.
\newblock {\em Journal of Computational Physics}, 225(2):1721 -- 1735, 2007.

\bibitem{lai_2000_ibm}
M.-C. Lai and C.~S. Peskin.
\newblock An immersed boundary method with formal second-order accuracy and
  reduced numerical viscosity.
\newblock {\em Journal of Computational Physics}, 160(2):705 -- 719, 2000.

\bibitem{liu_2003_phasefield}
C.~Liu and J.~Shen.
\newblock A phase field model for the mixture of two incompressible fluids and
  its approximation by a fourier-spectral method.
\newblock {\em Physica D: Nonlinear Phenomena}, 179(3–4):211 -- 228, 2003.

\bibitem{lombard_2015_sim}
J.-E.~W. Lombard, D.~Moxey, S.~J. Sherwin, J.~F.~A. Hoessler, S.~Dhandapani,
  and M.~J. Taylor.
\newblock Implicit large-eddy simulation of a wingtip vortex.
\newblock {\em AIAA Journal}, pages 1--13, Nov. 2015.

\bibitem{lynch_1964_tensors}
R.~Lynch, J.~Rice, and D.~Thomas.
\newblock Direct solution of partial difference equations by tensor product
  methods.
\newblock {\em Numerische Mathematik}, 6(1):185--199, 1964.

\bibitem{manna_2004_preconditioning}
M.~Manna, A.~Vacca, and M.~O. Deville.
\newblock Preconditioned spectral multi-domain discretization of the
  incompressible navier–stokes equations.
\newblock {\em Journal of Computational Physics}, 201(1):204 -- 223, 2004.

\bibitem{merzari_2013_sim}
E.~Merzari, W.~Pointer, and P.~Fischer.
\newblock Numerical simulation and proper orthogonal decomposition of the flow
  in a counter-flow t-junction.
\newblock {\em Journal of Fluids Engineering}, 135(9):091304, 2013.

\bibitem{patera_1984_sem}
A.~T. Patera.
\newblock A spectral element method for fluid dynamics: laminar flow in a
  channel expansion.
\newblock {\em Journal of Computational Physics}, 54(3):468 -- 488, 1984.

\bibitem{peskin_2002_ibm}
C.~S. Peskin.
\newblock The immersed boundary method.
\newblock {\em Acta Numerica}, 11:479--517, 1 2002.

\bibitem{pflaum_2008_multigrid}
C.~Pflaum.
\newblock A multigrid conjugate gradient method.
\newblock {\em {Applied Numerical Mathematics}}, {58}({12}):{1803--1817}, {DEC}
  {2008}.

\bibitem{shewchuk_1994_cg}
J.~R. Shewchuk.
\newblock An introduction to the conjugate gradient method without the
  agonizing pain.
\newblock Technical report, Pittsburgh, PA, USA, 1994.

\bibitem{trottenberg_2001_multigrid}
U.~Trottenberg, C.~Oosterlee, and A.~Sch{\"u}ller.
\newblock {\em Multigrid}.
\newblock Academic Press, 2001.

\bibitem{uhlmann_2005_ibm}
M.~Uhlmann.
\newblock An immersed boundary method with direct forcing for the simulation of
  particulate flows.
\newblock {\em Journal of Computational Physics}, 209(2):448 -- 476, 2005.

\bibitem{yakovlev_2015_hdg}
S.~Yakovlev, D.~Moxey, R.~Kirby, and S.~Sherwin.
\newblock To {CG} or to {HDG}: A comparative study in 3d.
\newblock {\em Journal of Scientific Computing}, pages 1--29, 2015.

\end{thebibliography}

\end{document}